\documentclass[12pt]{article}

\usepackage{amsmath}
\usepackage{amsfonts}
\usepackage{bigints}

\begin{document}

\renewcommand{\P}{\par \rm Proof:\ }
\newcommand{\Pe}{$\hfill{\Box}$\bigskip}
\newcommand{\Pes}{$\hfill{\Box}$}
\newcommand{\A}{\mbox{${{{\cal A}}}$}}
\newcommand{\Mm}{{\cal{M}}^{\mu}}
\newcommand{\Mn}{{\cal{M}}^{\nu}}


\author{Attila Losonczi}
\title{Measuring sets by means}

\date{16 August 2017}

\newtheorem{thm}{\qquad Theorem}[section]
\newtheorem{prp}[thm]{\qquad Proposition}
\newtheorem{lem}[thm]{\qquad Lemma}
\newtheorem{cor}[thm]{\qquad Corollary}
\newtheorem{rem}[thm]{\qquad Remark}
\newtheorem{ex}[thm]{\qquad Example}
\newtheorem{df}[thm]{\qquad Definition}
\newtheorem{prb}{\qquad Problem}

\maketitle

\begin{abstract}

\noindent

We are going to classify sets by a given mean in two ways. Firstly we study small and big sets regarding a given mean. Secondly we study sets that have the same weight according to a mean. We also generalize the notion of roundness and get another way to compare subsets by a mean.

\noindent
\footnotetext{\noindent
AMS (2010) Subject Classifications: 40A05, 40G05, 28A10, 28A78 \\

}

\end{abstract}

\section{Introduction}
This paper can be considered as a natural continuation of the investigations started in \cite{lamis} and \cite{lamisii} where we started to build the theory of means on infinite sets. An ordinary mean is for calculating the mean of two (or finitely many) numbers. This can be extended in many ways in order to get a more general concept where we have a mean on some infinite subsets of $\mathbb{R}$. The various general properties of such means, the relations among those means were studied thoroughly in \cite{lamis} and \cite{lamisii}. 

In this paper we make some efforts to qualitatively measure and compare some sets using a given mean. Our aim is to to classify the sets in the domain of a mean in two ways. Firstly we are going to identify small and big sets regarding a mean. Roughly speaking a set is considered small if adding to or subtracting from any other set it does not change the mean value of that set. On the contrary a set is big if for any set it changes its mean value when adding to or subtracting from it if it is positioned far enough from that set.

Secondly we are going to specify when two sets have the same weight regarding a given mean. Roughly it happens if the sets are moved far from each other then their mean starts to behave as if they were points. We define three such weight types.

In the last section we are dealing with a generalization of roundness. It turns out that it cannot be considered of a property of a mean, instead it can be seen as a property of the underlying set. Namely it says that the mean value of the set cuts the set into two equally weighted parts. 

\medskip

The investigation started in this paper can be applied to any type of mean on infinite sets however in most of the cases they esentially regard for arithmetic type means. In the sequel we will define that we require from "arithmetic typeness".

\subsection{Basic notions and notations}

For $K\subset\mathbb{R},\ y\in\mathbb{R}$ let us use the notation $K^{-y}=K\cap(-\infty,y],K^{+y}=K\cap[y,+\infty).$

If $H\subset\mathbb{R},x\in\mathbb{R}$ then set $H+x=\{h+x:h\in H\}$.

We use the convention that this operation $+$ has to be applied prior to the set theoretical operations, e.g. $H\cup K\cup L+x=H\cup K\cup (L+x)$.

$cl(H), H'$ will denote the closure and accumulation points of $H\subset\mathbb{R}$ respectively. Let $\varliminf H=\inf H',\ \varlimsup H=\sup H'$ for infinite bounded $H$. Let $H^{(1)}=H'$ and $H^{(n+1)}=(H^{(n)})'\ (n\geq1)$.

\medskip

Usually ${\cal{K}},{\cal{M}}$ will denote means, $Dom({\cal{K}})$ denotes the domain of ${\cal{K}}$. Please note that our means always deal with bounded sets only.

\medskip

Let us recall some definitions from \cite{lamis} and \cite{lamisii} that regards for means on infinite sets. 

A mean ${\cal{K}}$ is called \textbf{internal} if $\inf H\leq {\cal{K}}(H)\leq\sup H.$ It is \textbf{strongly internal} if $\varliminf H\leq {\cal{K}}(H)\leq\varlimsup H$ where $\varliminf H=\min H',\varlimsup H=\max H'$.

$\cal{K}$ is \textbf{monotone} if $\sup H_1\leq\inf H_2$ implies that ${\cal{K}}(H_1)\leq {\cal{K}}(H_1\cup H_2)\leq {\cal{K}}(H_2)$. 

${\cal{K}}$ is \textbf{strong monotone} if $\cal{K}$ is strong internal and $\varlimsup H_1\leq\varliminf H_2$ implies that ${\cal{K}}(H_1)\leq {\cal{K}}(H_1\cup H_2)\leq {\cal{K}}(H_2)$. 

$\cal{K}$ is \textbf{disjoint-monotone} if $H_1\cap H_2=\emptyset,{\cal{K}}(H_1)\leq{\cal{K}}(H_2)$ then ${\cal{K}}(H_1)\leq{\cal{K}}(H_1\cup H_2)\leq {\cal{K}}(H_2)$.

${\cal{K}}$ is \textbf{union-monotone} if $B\cap C=\emptyset$,
${\cal{K}}(A)\leq{\cal{K}}(A\cup B),{\cal{K}}(A)\leq{\cal{K}}(A\cup C)$ implies ${\cal{K}}(A)\leq{\cal{K}}(A\cup B\cup C)$ 
and 
${\cal{K}}(A\cup B)\leq{\cal{K}}(A),{\cal{K}}(A\cup C)\leq {\cal{K}}(A)$ implies ${\cal{K}}(A\cup B\cup C)\leq{\cal{K}}(A)$.
Moreover if any of the inequalities on the left hand side is strict then so is the inequality on the right hand side.

${\cal{K}}$ is \textbf{d-monotone} if $L,B\in Dom({\cal{K}}),L\cap B=(L\cup B)\cap (B+x)=\emptyset$ then
${\cal{K}}(L)<{\cal{K}}(L\cup B), x>0$ implies ${\cal{K}}(L\cup B)<{\cal{K}}(L\cup B\cup B+x)$
and

${\cal{K}}(L)>{\cal{K}}(L\cup B), x<0$ implies ${\cal{K}}(L\cup B)>{\cal{K}}(L\cup B\cup B+x)$.

 The mean is \textbf{shift invariant} if $x\in\mathbb{R}, H\in Dom({\cal{K}})$ then $H+x\in Dom({\cal{K}}),\ {\cal{K}}(H+x)={\cal{K}}(H)+x$.

Let $H\in Dom({\cal{K}}), x\in\mathbb{R}$ such that $\varlimsup H\leq\varliminf H+x$ or $\varlimsup H+x\leq\varliminf H$. Then ${\cal{K}}$ is called \textbf{self-shift-invariant} if ${\cal{K}}(H\cup (H+x))={\cal{K}}(H)+\frac{x}{2}$.

If $x\in\mathbb{R},H_1\cup H_2\in Dom({\cal{K}}),H_1\cup(H_2+x)\in Dom({\cal{K}}), H_1\cap H_2=H_1\cap H_2+x=\emptyset$ implies that $sign({\cal{K}}(H_1\cup(H_2+x))-{\cal{K}}(H_1\cup H_2))=sign(x)$ and $$|{\cal{K}}(H_1\cup(H_2+x))-{\cal{K}}(H_1\cup H_2)|\leq |x|$$ then ${\cal{K}}$ is called \textbf{part-shift-invariant}.

\medskip

Throughout this paper function $\A()$ will denote the arithmetic mean of any number of variables.

\begin{df}Let $H\subset\mathbb{R}$ be infinite and bounded.
$${\cal{M}}^{lis}(H)=\frac{\varliminf H+\varlimsup H}{2}$$
\end{df}

\begin{df}\label{davg}Let $\mu^s$ denote the s-dimensional Hausdorff measure ($0\leq s\leq 1$). If $0<\mu^s(H)<+\infty$ (i.e. $H$ is an $s$-set) and $H$ is $\mu^s$ measurable then $$Avg(H)=\frac{\int\limits_H x\ d\mu^s}{\mu^s(H)}.$$
If $0\leq s\leq 1$ then set $Avg^s=Avg|_{\{\text{measurable s-sets}\}}$. E.g. $Avg^1$ is $Avg$ on all Lebesgue measurable sets with positive measure.
\end{df}

\begin{df}If $cl(H-H')=H$ then let 
$${\cal{M}}^{iso}(H)=\lim_{\delta\to 0+0}\A(H-S(H',\delta))$$ 
if it exists.
\end{df}

\begin{df}Let $H\subset\mathbb{R}$. Let $lev(H)=n\in\mathbb{N}$ if $H^{(n+1)}=\emptyset$ and $H^{(n)}\ne\emptyset$. Otherwise let $lev(H)=+\infty$.
\end{df}

\begin{df}Let $H\subset\mathbb{R},\ lev(H)=n$. Let ${\cal{M}}^{acc}(H)=\A(H^{(n)})$. 
\end{df}

\begin{df}Let $H\in Dom({\cal{K}})$. Let
$\varliminf_{\cal{K}}H=\sup\{x:{\cal{K}}(H)={\cal{K}}(H^{+x})\}$, 
$\varlimsup_{\cal{K}}H=\inf\{x:{\cal{K}}(H)={\cal{K}}(H^{-x})\}$ be the liminf and limsup of $H$ with respect to the mean ${\cal{K}}$.
\end{df}

\smallskip

The means considered in this paper are more or less arithmetic type ones at least in the sense that we always require that they have properties shift-invariance and self-shift-invariance. However we will always mention if we use any of those properties.

\section{Small and big sets}
We are going to identify a kind of small sets regarding a given mean. 

\begin{df}Let ${\cal{K}}$ be a mean, $H\in Dom({\cal{K}})$. Let $S_{\cal{K}}(H)=\{V\subset\mathbb{R}$ bounded$:H\cup (V+x)\in Dom({\cal{K}}),H-(V+x)\in Dom({\cal{K}}),{\cal{K}}(H\cup V+x)={\cal{K}}(H-(V+x))={\cal{K}}(H)\ \forall x\in\mathbb{R}\}$.
\end{df}

We can consider $S_{\cal{K}}(H)$ as the small sets to $H$ regarding ${\cal{K}}$. Obviously $S_{\cal{K}}(H)$ is shift-invariant i.e. $V\in S_{\cal{K}}(H)$ iff $V+x\in S_{\cal{K}}(H)\ (x\in\mathbb{R})$.

\begin{df}Let ${\cal{K}}$ be a mean. Let $S_{\cal{K}}=\bigcap \{S_{\cal{K}}(H):H\in Dom({\cal{K}})\}$.
\end{df}

We can consider $S_{\cal{K}}$ as the small sets regarding ${\cal{K}}$. Clearly $S_{\cal{K}}$ is shift-invariant.

\begin{prp}If $Dom({\cal{K}})$ is shift-invariant then $S_{\cal{K}}\cap Dom({\cal{K}})=\emptyset$.
\end{prp}
\P Let $H\in S_{\cal{K}}\cap Dom({\cal{K}})$. Then $H+x\in S_{\cal{K}}\cap Dom({\cal{K}})\ (x\in\mathbb{R})$. Hence ${\cal{K}}(H)={\cal{K}}(H+x\cup H)={\cal{K}}(H+x)$. Then $H$ being bounded provides a contradiction.
\Pes

\begin{prp}If ${\cal{K}}$ is finite-independent then $S_{\cal{K}}$ contains all finite sets. More generally If ${\cal{K}}$ is ${\cal{I}}$-independent for a shift-invariant ideal ${\cal{I}}$ and $Dom({\cal{K}})\cap{\cal{I}}=\emptyset$ then ${\cal{I}}\subset S_{\cal{K}}$. \Pes
\end{prp}

\begin{ex}$S_{Avg^1}(H)=S_{Avg^1}=\{$the sets with 0 Lebesgue measure$\}\ (H\in Dom(Avg^1))$. \Pes
\end{ex}

\begin{ex}$S_{{\cal{M}}^{lis}}(H)=S_{{\cal{M}}^{lis}}=\{$finite sets$\}\ (H\in Dom({\cal{M}}^{lis}))$.
\end{ex}

\begin{lem}\label{lskcu}$S_{\cal{K}}$ is closed for union.
\end{lem}
\P If $V_1,V_2\in S_{\cal{K}}, x\in\mathbb{R}$ then ${\cal{K}}(H\cup (V_1\cup V_2)+x)={\cal{K}}(H\cup (V_1+x)\cup (V_2+x))={\cal{K}}(H\cup V_1+x)={\cal{K}}(H)$. Similarly ${\cal{K}}(H - (V_1\cup V_2)+x)={\cal{K}}(H - (V_1+x) - (V_2+x))={\cal{K}}(H - V_1+x)={\cal{K}}(H)$. 
This means that $V_1\cup V_2\in S_{\cal{K}}$.
\Pes

\begin{prp}If ${\cal{K}}$ is monotone, union-monoton, shift-invariant, part-shift-invariant then $S_{\cal{K}}$ is a shift-invariant ideal.
\end{prp}
\P By \ref{lskcu} we only have to show that $S_{\cal{K}}$ is descending.

Let $V_1\in S_{\cal{K}}$ and suppose that $V_2\subset V_1$ and $V_2\notin S_{\cal{K}}$. Let $V_1=V_2\cup^* V_3$. Assume that ${\cal{K}}(H)<{\cal{K}}(H\cup V_2)$. Find $x>0$ such that $\sup H-x\leq\inf V_3$. Then by monotonicity ${\cal{K}}(H-x)\leq{\cal{K}}(H-x\cup V_3)$. By part-shift-invariance we get ${\cal{K}}(H\cup V_2)-{\cal{K}}(H-x\cup V_2)\leq x$. Then ${\cal{K}}(H)-{\cal{K}}(H-x\cup V_2)< x$ which gives that ${\cal{K}}(H-x)<{\cal{K}}(H-x\cup V_2)$. Then by union-monotonicity we end up with ${\cal{K}}(H-x)<{\cal{K}}(H-x\cup V_2\cup V_3)={\cal{K}}(H-x\cup V_1)={\cal{K}}(H-x)$ which is a contradiction.

If ${\cal{K}}(H\cup V_2)<{\cal{K}}(H)$ then a very similar argument can be applied.
\Pe

We can now define disjointness regarding ${\cal{K}}$.

\begin{df}\label{dsbskd}Let $H_1,H_2\in Dom({\cal{K}})$. We say that $H_1,H_2$ are ${\cal{K}}$-disjoint if $H_1\cap H_2\in S_{\cal{K}}$. 
We call $H_1,H_2$ weak-${\cal{K}}$-disjoint if $H_1\cap H_2\in S_{\cal{K}}(H_1)\cap S_{\cal{K}}(H_2)$. 
\end{df}

Now let us define the big sets regarding a given mean.

\begin{df}Let ${\cal{K}}$ be a mean, $H\in Dom({\cal{K}})$. Let $B_{\cal{K}}(H)=\{V\in Dom({\cal{K}}):{\cal{K}}(V\cup(H+x))={\cal{K}}(V-(H+x))={\cal{K}}(V)\ \forall  x\in\mathbb{R}\}$.
\end{df}

We can consider $B_{\cal{K}}(H)$ as the big sets to $H$ regarding ${\cal{K}}$.

\begin{df}$H\in Dom({\cal{K}})$ is called a big set regarding ${\cal{K}}$ if $B_{\cal{K}}(H)=\emptyset$. Equivalently $\forall L\in Dom({\cal{K}})\ \exists x\in\mathbb{R}$ such that ${\cal{K}}(L\cup(H+x))\ne{\cal{K}}(L)$ or ${\cal{K}}(L)\ne{\cal{K}}(L-(H+x))$. Let us denote the set of big sets by $B_{\cal{K}}$.
\end{df}

One can readily check:

\begin{prp}\label{psmbieq}Let $H,V\in Dom({\cal{K}})$. Then $$V\in S_{\cal{K}}(H)\Longleftrightarrow H\in B_{\cal{K}}(V).$$
\end{prp}

\begin{ex}If $H$ is bounded infinite set then $B_{{\cal{M}}^{lis}}(H)=\emptyset$ since if $\varlimsup L<\varliminf H+x\leq\varlimsup H+x$ then ${\cal{K}}(L\cup(H+x))\ne{\cal{K}}(L)$. 
That implies that $B_{{\cal{M}}^{lis}}$ consists of all bounded infinite sets.
\end{ex}

\begin{ex}Similarly we get that $B_{Avg^1}=Dom(Avg^1)=\{$bounded sets with positive Lebesgue measure$\}$ because $B_{Avg^1}(H)=\emptyset\ \forall H\in Dom(Avg^1)$.
\end{ex}

\begin{ex}It can happen that there is no big set at all for a mean. Let ${\cal{K}}=Avg$ restricted to sets that are $s$-sets with $s<1$ (so simply leave out the 1-sets). Let $H\in Dom({\cal{K}})$ be an $s$-set ($s<1$). Then find a $V\in Dom({\cal{K}})$ that is an $s'$-set with $s<s'<1$. Clearly ${\cal{K}}(V\cup(H+x))={\cal{K}}(V-(H+x))={\cal{K}}(V)\ \forall  x\in\mathbb{R}$ i.e. $V\in B_{\cal{K}}(H)\ne\emptyset$.
\end{ex}

\begin{ex}Let ${\cal{K}}={\cal{M}}^{acc}$ and $H\in Dom({\cal{K}})$. Then $S_{\cal{K}}(H)=\{L\in Dom({\cal{K}}):lev(L)<lev(H)\}\cup\{$finite sets$\},B_{\cal{K}}(H)=\{L\in Dom({\cal{K}}):lev(L)>lev(H)\}$.
Those imply that $S_{\cal{K}}=\{$finite sets$\},\ B_{\cal{K}}=\emptyset$. 
\end{ex}

\begin{prp}\label{pskbkiso}Let ${\cal{K}}={\cal{M}}^{iso},\ H_1,H_2\in Dom({\cal{K}})$. Set $n_{\epsilon}=|H_1-S(H_1',\epsilon)|,m_{\epsilon}=|H_2-S(H_2',\epsilon)|\ (\epsilon>0)$.

Then $H_1\in S_{\cal{K}}(H_2)$ iff $\lim_{\epsilon\to0+0}\frac{n_{\epsilon}}{m_{\epsilon}}=0$. Similarly $H_1\in B_{\cal{K}}(H_2)$ iff $\lim_{\epsilon\to0+0}\frac{n_{\epsilon}}{m_{\epsilon}}=+\infty$.
\end{prp}
\P Suppose first that $H_1\in S_{\cal{K}}(H_2)$. We can assume that $1+\sup H_2<\inf H_1$. Set $H=H_1\cup H_2$. Let $1>\epsilon>0$. Then 
\begin{equation}\label{eqiso}
\A(H-S(H',\epsilon))=\frac{n_{\epsilon}\A(H_1-S(H_1',\epsilon))+m_{\epsilon}\A(H_2-S(H_2',\epsilon))}{n_{\epsilon}+m_{\epsilon}}
\end{equation}
that is equivalent to
$$\frac{n_{\epsilon}}{m_{\epsilon}}\Big(\A(H-S(H',\epsilon)) - \A(H_1-S(H_1',\epsilon))\Big)=\A(H_2-S(H_2',\epsilon)) - \A(H-S(H',\epsilon)).$$
By assumption $\A(H-S(H',\epsilon))\to {\cal{K}}(H_2)$ when $\epsilon\to0+0$. It means that the right hand side tends to 0, while $\A(H-S(H',\epsilon)) - \A(H_1-S(H_1',\epsilon))$ tends to a non zero number (${\cal{K}}(H_2)-{\cal{K}}(H_1)>1$). It gives that $\frac{n_{\epsilon}}{m_{\epsilon}}$ must tend to $0$ as well.

Now suppose that $\lim_{\epsilon\to0+0}\frac{n_{\epsilon}}{m_{\epsilon}}=0$. By (\ref{eqiso}) we get that
$$\A(H-S(H',\epsilon))=\frac{\frac{n_{\epsilon}}{m_{\epsilon}}\A(H_1-S(H_1',\epsilon))+\A(H_2-S(H_2',\epsilon))}{\frac{n_{\epsilon}}{m_{\epsilon}}+1}$$
which yields that $\A(H-S(H',\epsilon))\to  {\cal{K}}(H_2)$ when $\epsilon\to0+0$.

The "big" case is similar or can be referred to \ref{psmbieq}.
\Pes

\begin{prp}Let ${\cal{K}}={\cal{M}}^{iso}$. Then for all $H_2\in Dom({\cal{K}})\ S_{\cal{K}}(H_2)\ne\emptyset$ and $B_{\cal{K}}(H_2)\ne\emptyset$.
\end{prp}
\P Let us use notations of \ref{pskbkiso}. 

Let $H_2\in Dom({\cal{K}})$ be arbitrary. It is easy to construct an $H_1\in B_{\cal{K}}(H_2)$ in the following way. 
Create $2|H_2-S(H_2',\frac{1}{2})|$ many points in $S(H_2',\frac{1}{1})-S(H_2',\frac{1}{2})$. Then create $3|H_2-S(H_2',\frac{1}{3})|$ many points in $S(H_2',\frac{1}{2})-S(H_2',\frac{1}{3})$. Generally create $n|H_2-S(H_2',\frac{1}{n})|$ many points in $S(H_2',\frac{1}{n-1})-S(H_2',\frac{1}{n})$. Let $H_1$ consist of all those points.

Let us show that for any set $H_2 \in Dom({\cal{K}})$ we can construct an $H_1\in S_{\cal{K}}(H_2)$. Take a point in $S(H_2',\frac{1}{1})-S(H_2',\frac{1}{2})$. Then find $k_1\in\mathbb{N}$ such that $\frac{1}{m_{\frac{1}{k_1}}}<\frac{1}{2}$. Then take another single point in $S(H_2',\frac{1}{k_1})-S(H_2',\frac{1}{k_1+1})$. Then find $k_2\in\mathbb{N}$ such that $\frac{2}{m_{\frac{1}{k_2}}}<\frac{1}{3}$. Then take another single point in $S(H_2',\frac{1}{k_2})-S(H_2',\frac{1}{k_2+1})$. Then find $k_3\in\mathbb{N}$ such that $\frac{3}{m_{\frac{1}{k_3}}}<\frac{1}{4}$. And so on. Let $H_1$ consist of all those points.
\Pes

\begin{cor}$S_{{\cal{M}}^{iso}}=\{$finite sets$\},\ B_{{\cal{M}}^{iso}}=\emptyset$. 
\end{cor}
\P Only the first statement needs proof. Let $H$ be a bounded infinite set. Take one of its accumulation point, say $a\in H$ and an $(a_n)$ such that $a_n\to a$ and $a_n\in H\ (n\in\mathbb{N})$. Let $K=\{a_n:n\in\mathbb{N}\}$. Then $H\not\in S_{{\cal{M}}^{iso}}(K)$ hence $H\not\in S_{{\cal{M}}^{iso}}$.
\Pes

\begin{prp}If $S_{\cal{K}}$ is descending, $B\in B_{\cal{K}}, S\in S_{\cal{K}}$ then $B-S\in B_{\cal{K}}$.
\end{prp}
\P Let $L\in Dom({\cal{K}}),x\in\mathbb{R}$ such that ${\cal{K}}(L\cup(B+x))\ne{\cal{K}}(L)$ or ${\cal{K}}(L-(B+x))\ne {\cal{K}}(L)$.

${\cal{K}}(L\cup((B-S)+x))={\cal{K}}(L\cup(B+x)-((S+x)-L))={\cal{K}}(L\cup(B+x))$ because $(S+x)-L\in S_{\cal{K}}$.

${\cal{K}}(L-((B-S)+x))={\cal{K}}((L-(B+x))\cup(L\cap S+x))={\cal{K}}(L-(B+x))$ because $L\cap S+x\in S_{\cal{K}}$.

Hence either ${\cal{K}}(L\cup((B-S)+x))\ne{\cal{K}}(L)$ or ${\cal{K}}(L-((B-S)+x))\ne{\cal{K}}(L)$.
\Pes

\begin{lem}\label{prpunm}If ${\cal{K}}$ is monotone, d-monotone, union-monotone then $B\in B_{\cal{K}}$ iff $\forall L\in Dom({\cal{K}})\ \exists x\in\mathbb{R}$ such that ${\cal{K}}(L\cup(B+x))\ne{\cal{K}}(L)$.
\end{lem}
\P We have to prove that 

$L\in Dom({\cal{K}}),{\cal{K}}(L)\ne{\cal{K}}(L-(B+x)) \Longrightarrow $
$\exists y\in\mathbb{R}$ such that ${\cal{K}}(L\cup(B+y))\ne{\cal{K}}(L)$.

Suppose ${\cal{K}}(L-(B+x))<{\cal{K}}(L)$. Let $B_1=L\cap (B+x), L_1=L-B_1$. Then $L_1=L-(B+x),L=L_1\cup B_1$ hence ${\cal{K}}(L_1)<{\cal{K}}(L_1\cup B_1)$. If $y>0,(L_1\cup B_1)\cap B_1+y=\emptyset$ then ${\cal{K}}(L)={\cal{K}}(L_1\cup B_1)<{\cal{K}}(L_1\cup B_1\cup B_1+y)={\cal{K}}(L\cup B_1+y)$ by d-monotonicity. 

Choose $y$ such that $\sup L\leq\inf (B-B_1)+y$. Then ${\cal{K}}(L)\leq{\cal{K}}(L\cup (B-B_1)+y)$ by monotonicity. Therefore by union-monotonicity ${\cal{K}}(L)<{\cal{K}}(L\cup (B_1+y)\cup (B-B_1)+y))={\cal{K}}(L\cup (B+y))$.

The case when ${\cal{K}}(L)<{\cal{K}}(L-(B+x))$ can be handled similarly.
\Pes

\begin{lem}\label{lbxy}Let ${\cal{K}}$ is monotone, union-monotone, part-shift-monotone. Let ${\cal{K}}(L)<{\cal{K}}(L\cup(B_1+x))$. Then there is $z\in\mathbb{R}$ such that $y\geq z$ implies that ${\cal{K}}(L)<{\cal{K}}(L\cup(B_1+y))$.
\end{lem}
\P Let $y>x$ such that $\sup L<\inf B_1+y$. 

Let $B_3=(B_1+x)-L,B_4=(B_1+x)\cap L$. Then $B_3\cap L=B_3+y-x\cap L=\emptyset$ hence by part-shift-monotonicity we get that ${\cal{K}}(L)<{\cal{K}}(L\cup B_1+x)={\cal{K}}(L\cup B_3)\leq{\cal{K}}(L\cup (B_3+y-x))$. By monotonicity ${\cal{K}}(L)\leq{\cal{K}}(L\cup B_4+y-x)$. Finally by union-monotonicity ${\cal{K}}(L)<{\cal{K}}(L\cup(B_3+y-x)\cup(B_4+y-x))={\cal{K}}(L\cup B_1+y)$.
\Pes

Obviously we can formulate a similar lemma for the opposite inequality.

\begin{thm}If ${\cal{K}}$ is monotone, d-monotone, union-monotone, part-shift-monotone then $B_1\in B_{\cal{K}},B_1\subset B_2$ implies that $B_2\in B_{\cal{K}}$ i.e. $B_{\cal{K}}$ is ascending.
\end{thm}
\P By Lemma \ref{prpunm} we have to handle the union part only.

Let $L\in Dom({\cal{K}})$. Choose $x\in\mathbb{R}$ such that ${\cal{K}}(L\cup(B_1+x))\ne{\cal{K}}(L)$. Suppose ${\cal{K}}(L)<{\cal{K}}(L\cup(B_1+x))$ (the other case can be proved similarly).

Find $y>x$ such that $\sup L\leq \inf B_2+y$ and $\sup L<\inf B_1+y$. Then ${\cal{K}}(L)\leq{\cal{K}}(L\cup (B_2-B_1)+y)$ by monotonicity. By the proof of Lemma \ref{lbxy} ${\cal{K}}(L)<{\cal{K}}(L\cup B_1+y)$.

By union monotonicity ${\cal{K}}(L)<{\cal{K}}(L\cup(B_1+y)\cup(B_2-B_1)+y)={\cal{K}}(L\cup(B_2+y))$.
\Pe

Now we introduce the notion of no-small and no-big.

\begin{df}Let $H,V\in Dom({\cal{K}})$. We say that $H,V$ are comparable regarding ${\cal{K}}$ whenever $V\not\in S_{\cal{K}}(H)\cup B_{\cal{K}}(H)$.
\end{df}

\begin{rem}This relation is symmetric i.e. the condition is equivalent to $H\not\in S_{\cal{K}}(V)\cup B_{\cal{K}}(V)$.
\end{rem}
\P \ref{psmbieq}.
\Pes

\begin{ex}For ${\cal{K}}={\cal{M}}^{acc}\ H,V\in Dom({\cal{K}})$ are comparable iff $lev(H)=lev(V)$. 
\end{ex}

\begin{ex}For ${\cal{K}}=Avg^1$ or ${\cal{K}}={\cal{M}}^{lis}$ all $H,V\in Dom({\cal{K}})$ are comparable.
\end{ex}

\section{Sets of equal weight}
We are now going to measure sets by a mean ${\cal{K}}$ at least in the sense that we could say that somehow they have equal weight by ${\cal{K}}$. Our first guess could be that ${\cal{K}}(H_1\cup^* H_2)=\frac{{\cal{K}}(H_1)+{\cal{K}}(H_2)}{2}$ would work as a criteria for that. However a simple example shows that does not.

\begin{ex}Let $H_1=\{1,2\}, H_2=\{0.5,1,3\}$. Then $\A(H_1\cup^* H_2)=\frac{\A(H_1)+\A(H_2)}{2}$ holds.
\end{ex}

The following proposition shows that e.g. for ${\cal{K}}=\A$ or $Avg$ this condition cannot be used. 

\begin{prp}Let ${\cal{K}}$ be disjoint-monotone, $H_1,H_2\in Dom({\cal{K}}),\ H_1\cap H_2=\emptyset$ and  ${\cal{K}}(H_1)={\cal{K}}(H_2)$. Then ${\cal{K}}(H_1\cup^* H_2)=\frac{{\cal{K}}(H_1)+{\cal{K}}(H_2)}{2}$. 
\Pes
\end{prp}

Hence we need something stronger. It is based on the observation that such "equal" sets behave like single points when they are moved very far from each other. 

\begin{df}\label{dsewm}Let ${\cal{K}}$ be a shift-invariant, monotone mean, $H_1,H_2\in Dom({\cal{K}})$. We say that $H_1,H_2$ have \textbf{equal weight in bound} regarding ${\cal{K}}$ if  
$$\sup_{x\in\mathbb{R}}\left\lvert{\cal{K}}(H_1\cup (H_2+x))-\frac{{\cal{K}}(H_1)+{\cal{K}}(H_2+x)}{2}\right\lvert<+\infty.$$

We say that $H_1,H_2$ have \textbf{equal weight in limit} regarding ${\cal{K}}$ if 
$$\lim_{x\to\infty}\left\lvert{\cal{K}}(H_1\cup (H_2+x))-\frac{{\cal{K}}(H_1)+{\cal{K}}(H_2+x)}{2}\right\lvert=0.$$

We say that $H_1,H_2$ have \textbf{equal weight in equality} regarding ${\cal{K}}$ if $\varlimsup_{\cal{K}} H_1<\varliminf_{\cal{K}} H_2+x$ or $\varlimsup_{\cal{K}} H_2+x<\varliminf_{\cal{K}} H_1$ implies that
$${\cal{K}}(H_1\cup (H_2+x))=\frac{{\cal{K}}(H_1)+{\cal{K}}(H_2+x)}{2}.$$
\end{df}

\begin{prp}Clearly: $H_1,H_2$ have equal weight in equality regarding ${\cal{K}}\ \Longrightarrow\ H_1,H_2$ have equal weight in limit $\Longrightarrow\ H_1,H_2$ have equal weight in bound. \Pes
\end{prp}

\begin{prp}If ${\cal{K}}$ is shift-invariant then if $H_1,H_2$ have equal weight in bound (limit/equality) regarding ${\cal{K}}$ then so do $H_2,H_1$ i.e. this relation is symmetric.
\end{prp}
\P Observe that ${\cal{K}}(H_2\cup (H_1+x))={\cal{K}}(H_1\cup (H_2-x))+x$ and 
$$\frac{{\cal{K}}(H_2)+{\cal{K}}(H_1+x)}{2}=\frac{{\cal{K}}(H_1)+{\cal{K}}(H_2-x)}{2}+x.$$
Those imply that 
$${\cal{K}}(H_2\cup (H_1+x))-\frac{{\cal{K}}(H_2)+{\cal{K}}(H_1+x)}{2}={\cal{K}}(H_1\cup (H_2-x))-\frac{{\cal{K}}(H_1)+{\cal{K}}(H_2-x)}{2}$$
which gives the statement in bound and limit.

For "equality" let us note that $\varlimsup_{\cal{K}} H_2<\varliminf_{\cal{K}} H_1+x$ is equivalent to $\varlimsup_{\cal{K}} H_2-x<\varliminf_{\cal{K}} H_1$ (and similarly the other inequality).
\Pes

\begin{prp}If ${\cal{K}}$ is self-shift-invariant then $H$ has equal weight in equality to itself regarding ${\cal{K}}$ i.e. this relation is reflexive. \Pes
\end{prp}

\begin{prp}\label{psew1}If ${\cal{K}}=\A,\ H_1,H_2\in Dom({\cal{K}})$ are finite sets and $H_1,H_2$ have equal weight in bound regarding ${\cal{K}}$ then $|H_1|=|H_2|$.
\end{prp}
\P Let $H_1=\{h_1,\dots,h_n\},H_2=\{k_1,\dots,k_m\}$. Then ${\cal{K}}(H_1\cup (H_2+x))=\frac{\sum_1^n h_i+\sum_1^m k_j+mx}{n+m}$ and $\frac{{\cal{K}}(H_1)+{\cal{K}}(H_2+x)}{2}=\frac{\frac{\sum_1^n h_i}{n}+\frac{\sum_1^m k_j}{m}+x}{2}$. Set $C=\frac{\sum_1^n h_i+\sum_1^m k_j}{n+m}, D=\frac{\frac{\sum_1^n h_i}{n}+\frac{\sum_1^m k_j}{m}}{2}$. Then we get $|C+\frac{m}{n+m}x-(D+\frac{x}{2})|<K$ for a suitable $K$. Equivalenty $C-D-(\frac{m}{n+m}-\frac{1}{2})x$ is bounded that implies that $n=m$.
\Pes

\begin{prp}If ${\cal{K}}=Avg,\ H_1,H_2\in Dom({\cal{K}}), H_1$ is $s_1$-set, $H_2$ is $s_2$-set and $H_1,H_2$ have equal weight in bound regarding $Avg$ then $s_1=s_2,\ \mu^{s_1}(H_1)=\mu^{s_2}(H_2)$.
\end{prp}
\P If e.g. $s_1<s_2$ then ${\cal{K}}(H_1\cup (H_2+x))={\cal{K}}(H_2+x)={\cal{K}}(H_2)+x$ and then we got ${\cal{K}}(H_2)+x-(\frac{{\cal{K}}(H_1)+{\cal{K}}(H_2)}{2}+\frac{x}{2})$ is bounded that is impossible hence $s_1=s_2$.

If $s_1=0$ then we get the statement by Proposition \ref{psew1}. Let us assume that $s_1>0$. Let $\mu=\mu^{s_1}$. Then
$$Avg(H_1\cup (H_2+x))=\frac{\int_{H_1\cup (H_2+x)}y\ d\mu(y)}{\mu(H_1)+\mu(H_2)}=\frac{\mu(H_1)Avg(H_1)+\mu(H_2)Avg(H_2)+\mu(H_2)x}{\mu(H_1)+\mu(H_2)}.$$

$$\frac{Avg(H_1)+Avg(H_2+x)}{2}=\frac{Avg(H_1)+Avg(H_2)}{2}+\frac{x}{2}.$$

As in Proposition \ref{psew1} we get that $(\frac{1}{2}-\frac{\mu(H_2)}{\mu(H_1)+\mu(H_2)})x$ is bounded that immediately gives that $\mu^{s_1}(H_1)=\mu^{s_2}(H_2)$.
\Pes

\begin{rem}The proof also shows that for $Avg$ sets of equal weight in bound are also of equal weight in equality.
\end{rem}

\begin{prp}If ${\cal{K}}={\cal{M}}^{acc},\ H_1,H_2\in Dom({\cal{K}}), lev(H_1)=l_1, lev(H_2)=l_2$ and $H_1,H_2$ have equal weight in bound regarding ${\cal{M}}^{acc}$ then $l_1=l_2,\ |H_1^{(l_1)}|=|H_2^{(l_2)}|$.
\end{prp}
\P  If e.g. $l_1<l_2$ then ${\cal{K}}(H_1\cup (H_2+x))=\A(H_2^{(l_2)})={\cal{K}}(H_2+x)={\cal{K}}(H_2)+x$ and then we got ${\cal{K}}(H_2)+x-(\frac{{\cal{K}}(H_1)+{\cal{K}}(H_2)}{2}+\frac{x}{2})$ is bounded that is impossible hence $l_1=l_2$.

Then by Proposition \ref{psew1} we get $|H_1^{(l_1)}|=|H_2^{(l_2)}|$.
\Pes

\begin{prp}If ${\cal{K}}={\cal{M}}^{iso},\ H_1,H_2\in Dom({\cal{K}})$ and $H_1,H_2$ have equal weight in bound regarding ${\cal{M}}^{iso}$ then  $\lim_{\epsilon\to0+0}\frac{m_{\epsilon}}{n_{\epsilon}}=1$ where $n_{\epsilon}=|H_1-S(H_1',\epsilon)|,m_{\epsilon}=|H_2-S(H_2',\epsilon)|$.
\end{prp}
\P Let us take a sequence $(\epsilon_k)$ such that $\epsilon_k>0,\epsilon_k\to 0$ and $\frac{m_{\epsilon_k}}{n_{\epsilon_k}}\to l<+\infty$. Let $a=inf H_1,b=inf H_2,\ K\in\mathbb{R}$ such that $\sup H_1<a+K,\sup H_2<b+K$. Let ${\cal{K}}(H_1)=k_1=a+w,{\cal{K}}(H_2)=k_2=b+z$.

Let us denote the elements of $H_1-S(H_1',\epsilon)$ by $a_1,a_2,\dots$ and the elements of $H_2-S(H_2',\epsilon)$ by $b_1,b_2,\dots$. Then ${\cal{K}}(H_1\cup (H_2+x))=\lim_{\epsilon\to0+0}\frac{\sum_{i=1}^{n_\epsilon}a_i+\sum_{j=1}^{m_\epsilon}b_j+m_\epsilon x}{n_\epsilon+m_\epsilon}=v$.
We get that 
$$\frac{{n_{\epsilon_k}}a+{m_{\epsilon_k}}(b+x)}{n_{\epsilon_k}+m_{\epsilon_k}}=\frac{a+\frac{m_{\epsilon_k}}{n_{\epsilon_k}}(b+x)}{1+\frac{m_{\epsilon_k}}{n_{\epsilon_k}}}\leq\frac{\sum_{i=1}^{n_{\epsilon_k}}a_i+\sum_{j=1}^{m_{\epsilon_k}}b_j+m_{\epsilon_k} x}{n_{\epsilon_k}+m_{\epsilon_k}}\leq \frac{a+K+\frac{m_{\epsilon_k}}{n_{\epsilon_k}}(b+K+x)}{1+\frac{m_{\epsilon_k}}{n_{\epsilon_k}}}.$$

If $k\to\infty$ then 
$$\frac{a+l(b+x)}{1+l}\leq v\leq\frac{a+l(b+x)+Kl+K}{1+l}=\frac{a+l(b+x)}{1+l}+K$$

By assumption $v-\frac{k_1+k_2+x}{2}=v-\frac{a+b+x}{2}+\frac{w+z}{2}$ is bounded. Therefore $\frac{a+l(b+x)}{1+l}-\frac{a+b+x}{2}=\frac{(a-(b+x))(1-l)}{2(1+l)}$ has to be bounded as well. This yields that $l=1$.

If  we assumed that $\frac{m_{\epsilon_k}}{n_{\epsilon_k}}\to +\infty$ then take $\frac{n_{\epsilon_k}}{m_{\epsilon_k}}\to p=0$ instead and similarly one can show that this leads to a contradiction as $p=1$.
\Pes

\begin{rem}Actually we get that sets of equal weight in bound are also of equal weight in limit/equality.
\end{rem}
\P By the notation of the previous Proposition we get that
$$\frac{\sum_{i=1}^{n_{\epsilon}}a_i+\sum_{j=1}^{m_{\epsilon}}b_j+m_{\epsilon} x}{n_{\epsilon}+m_{\epsilon}}
=
\frac{\frac{\sum_{i=1}^{n_{\epsilon}}a_i}{n_{\epsilon}}+\frac{m_{\epsilon}}{n_{\epsilon}}\frac{\sum_{j=1}^{m_{\epsilon}}b_j}{m_{\epsilon}}+ x}{1+\frac{m_{\epsilon}}{n_{\epsilon}}}\to \frac{k_1+k_2+x}{2}.$$
\Pes

\begin{prp}If ${\cal{K}}={\cal{M}}^{lis},\ H_1,H_2\in Dom({\cal{K}})$ then $H_1,H_2$ always have equal weight in bound.

If $H_1,H_2$ are infinite, bounded and have equal weight in limit/equality regarding ${\cal{M}}^{lis}$ then $\varlimsup H_1-\varliminf H_1=\varlimsup H_2-\varliminf H_2$ .
\end{prp}
\P Let $H_1,H_2$ be infinite and let $\varliminf H_1=a,\varlimsup H_1=b,\varliminf H_2=c,\varlimsup H_2=d$. Then we get that $\frac{a+d+x}{2}-\frac{a+b+c+d+2x}{4}$ has to be bounded and that always holds.

Regarding limit/equality $\frac{a+d}{2}-\frac{a+b+c+d}{4}=0$ is equivalent to $d-c=b-a$.
\Pe

Towards transitivity of this "equal-weight" relation we can prove the following statetement that shows the required condition in raw form.

\begin{prp}Let $H_1,H_2$ and $H_2,H_3$ and $H_3,H_1$ are all sets of pair of equal weight in the same manner (in bound or limit or equality). Then 
$$\left\lvert{\cal{K}}(H_1\cup H_2+x)+{\cal{K}}(H_2+x\cup H_3+2x) - ({\cal{K}}(H_1\cup H_3+2x)+{\cal{K}}(H_2+x))\right\lvert\text{ is bounded}$$
$$\text{or }\lim_{x\to\infty}\left\lvert{\cal{K}}(H_1\cup H_2+x)+{\cal{K}}(H_2+x\cup H_3+2x) - ({\cal{K}}(H_1\cup H_3+2x)+{\cal{K}}(H_2+x))\right\lvert=0$$
$$\text{or }{\cal{K}}(H_1\cup H_2+x)+{\cal{K}}(H_2+x\cup H_3+2x) = {\cal{K}}(H_1\cup H_3+2x)+{\cal{K}}(H_2+x)$$
 holds respectively to the equal-weight type.
\end{prp}
\P We show it for equality type. ${\cal{K}}(H_1\cup H_2+x)+{\cal{K}}(H_2+x\cup H_3+2x) = \frac{{\cal{K}}(H_1)+{\cal{K}}(H_2+x)}{2} + \frac{{\cal{K}}(H_2+x)+{\cal{K}}(H_3+2x)}{2} = \frac{{\cal{K}}(H_1)+{\cal{K}}(H_3+2x)}{2} +{\cal{K}}(H_2+x)= {\cal{K}}(H_1\cup H_3+2x)+{\cal{K}}(H_2+x)$.
\Pe

We provide a sufficient condition for transitivity of the "equal-weight" relation. We prove it for the "limit" version, the others can be handled similarly.

\begin{prp}Let ${\cal{K}}$ be shift-invariant and self-shift-invariant. Let us use the notation $A\sim B$ if $A,B$ have equal weight in limit regarding ${\cal{K}}$.

Let the following hold: If $A,B,C,D\in Dom({\cal{K}})$ are pairwise disjoint sets and $A\sim B$ and $C\sim D$ then $A\cup C\sim B\cup D$ as well. 

In this case the "equal-weight in limit" relation is transitive.
\end{prp}
\P First observe that if $A\sim B$ then $A\sim B+x$ also holds.

Suppose that $A\sim B, B\sim C$ hold. Then for $x$ big enough $A\sim B+x, B+2x\sim C+3x$ also hold and the sets $A,B+x,B+2x,C+3x$ are pairwise disjoint. 
By asssumption $A\cup B+2x\sim B+x\cup C+3x$ and $A\cup C+3x\sim B+x\cup B+2x$. Which means that
$$\lim_{x\to\infty}\left\lvert{\cal{K}}(A\cup B+2x\cup (B+x\cup C+3x))-\frac{{\cal{K}}(A\cup B+2x)+{\cal{K}}(B+x\cup C+3x)}{2}\right\lvert=$$
$$\lim_{x\to\infty}\left\lvert\frac{{\cal{K}}(A\cup B+2x)+{\cal{K}}(B+x\cup C+3x)}{2}-\frac{\frac{{\cal{K}}(A)+{\cal{K}}(B+2x)}{2}+\frac{{\cal{K}}(B+x)+ {\cal{K}}(C+3x)}{2}}{2}\right\lvert=0.$$
Similarly
$$\lim_{x\to\infty}\left\lvert{\cal{K}}(A\cup C+3x\cup (B+x\cup B+2x))-\frac{{\cal{K}}(A\cup C+3x)+{\cal{K}}(B+x\cup B+2x)}{2}\right\lvert=$$
$$\lim_{x\to\infty}\left\lvert\frac{{\cal{K}}(A\cup C+3x)+{\cal{K}}(B+x\cup B+2x)}{2}-\frac{{\cal{K}}(A\cup C+3x)+\frac{{\cal{K}}(B+x)+ {\cal{K}}(B+2x)}{2}}{2}\right\lvert=0$$ where we used that $B+x\sim B+2x$ by self-shift-invariance.

These yield that 
$$\lim_{x\to\infty}\left\lvert{\cal{K}}(A\cup C+3x)-\frac{{\cal{K}}(A)+ {\cal{K}}(C+3x)}{2}\right\lvert=0$$
which is exactly that we had to prove.
\Pes

\subsection{On a roundness type notion}

Let us recall the notion of roundness.
\begin{df}A 2-variable mean $\circ$ is called round if it fulfils functional equation $(a\circ k)\circ (k\circ b)=k$ for $\forall a,b$ where $k=a\circ b$.
\end{df}

If we wanted to generalizise this notion for means acting on infinite sets then we would end up with something like this.

${\cal{K}}$ is round if the following holds. Let $k={\cal{K}}(H), k_1={\cal{K}}(H^{-k}), k_2={\cal{K}}(H^{+k})$. Then $k=\frac{k_1+k_2}{2}$.

It is easy to see that none of our means satisfy this too strong condition. E.g.

\begin{ex}$Avg$ is not round. 

Let $H=[0,2]\cup[4,5]$. Then $Avg(H)=\frac{2^2-0^2+5^2-4^2}{2\cdot 3}=\frac{13}{6},Avg([0,2])=1,Avg([4,5])=4.5$. Hence $\frac{13}{6}\ne\frac{1+4.5}{2}$.
\end{ex}

Actually roundness is not a property of a mean, instead it says something about the set. Roughly speaking it states that $H^{-k}$ and $H^{+k}$ have the same weight in some sense.
Hence it is better to say that $H\in Dom({\cal{K}})$ is round regarding ${\cal{K}}$ if the above property holds (for $H$ and ${\cal{K}}$).

\begin{df}$H\in Dom({\cal{K}})$ is round regarding ${\cal{K}}$ if $H^{-k},H^{+k}\in Dom({\cal{K}})$ and $\frac{k_1+k_2}{2}=k$ where $k={\cal{K}}(H), k_1={\cal{K}}(H^{-k}), k_2={\cal{K}}(H^{+k})$.
\end{df}

In the sequel we investigate how it behaves for some of the usual means.

\begin{prp}\label{pravg}Let $H$ be an s-set ($0< s\leq 1$) with $0<\mu^s(H)<+\infty$. Then $H$ is round regarding $Avg$ iff $\mu^s(H^{-k})=\mu^s(H^{+k})$ where $k=Avg(H)$.
\end{prp}
\P Let $\mu=\mu^s$. Then
$$Avg(H)=\frac{\mu(H^{-k})Avg(H^{-k})+\mu(H^{+k})Avg(H^{+k})}{\mu(H)}=\frac{Avg(H^{-k})+Avg(H^{+k})}{2}$$
which yields
$Avg(H^{-k})(2\mu(H^{-k})-\mu(H))+Avg(H^{+k})(2\mu(H^{+k})-\mu(H))=0$ that is equivalent to
$(Avg(H^{-k})-Avg(H^{+k}))(2\mu(H^{-k})-\mu(H))=0$. By strict strong internality of $Avg$ we know that $Avg(H^{-k})<k<Avg(H^{+k})$. Therefore $\mu(H^{-k})=\frac{\mu(H)}{2}$. 

Of course each transformation that we did can be reversed hence we proved the equivality of the two statements.
\Pes

\begin{prp}\label{paround}Let $H$ be finite. Then $H$ is round regarding $\A$ iff $|H^{-k}|=|H^{+k}|$ where $k=\A(H)$.
\end{prp}
\P Let $k=\A(H)$. First suppose that $k\not\in H$. Let $A=\A(\{h\in H:h<k\}),B=\A(\{h\in H:h>k\})$ and $m_1=|H^{-k}|,m_2=|H^{+k}|$. Then $H$ is round regarding $\A$ iff $\frac{A+B}{2}=\A(H)=\frac{m_1 A+m_2 B}{m_1+m_2}$ that is $m_2 A+m_1 B=m_1 A+m_2 B$ which gives that $(m_1-m_2)(A-B)=0$. But obviously $A<k<B$ hence $m_1=m_2$ that we had to prove.

If $k\in H$ then clearly $\A(H_1)=k$ where $H_1=H-\{k\}$. Let $m_1=|H_1^{-k}|,m_2=|H_1^{+k}|, a=\A(H_1^{-k}), b=\A(H_1^{+k})$. 

Firstly assume that $H$ is round.
Our aim is to show that $\frac{a+b}{2}=k$. That implies that $H_1$ is round hence by the previous argument we get that $|H_1^{-k}|=|H_1^{+k}|$ that is equivalent to $|H^{-k}|=|H^{+k}|$.

If $H$ is round then 
$$\frac{\frac{m_1 a+k}{m_1+1}+\frac{m_2 b+k}{m_2+1}}{2}=k \text{ where } k=\frac{m_1 a+k+m_2 b}{m_1+m_2+1}.$$
From the first we get $m_1 a+k+m_1 m_2 a+k m_2+m_2 b+k+m_1 m_2 b+k m_1=2k m_1 m_2+2k m_1+2k m_2+2k$ that is 
$$m_1 a+m_1 m_2 a+m_2 b+m_1 m_2 b=2k m_1 m_2+k m_1+k m_2.$$
From the second one we get 
$$m_1 a+m_2 b=k m_1+k m_2.$$
Subtracting the second from the first yields
$$m_1 m_2 a+m_1 m_2 b=2k m_1 m_2$$
that is $\frac{a+b}{2}=k.$

\medskip

Let us assume now that $|H^{-k}|=|H^{+k}|$. Then $m_1=m_2=m$ and $H_1$ is round which gives that $\frac{a+b}{2}=k.$
We want to prove that $H$ is round that means  
$$\frac{\frac{m a+k}{m+1}+\frac{m b+k}{m+1}}{2}=k$$
which is $ma+mb+2k=2mk+2k$ that is true.
\Pes

\begin{rem}As $Avg^0=\A$ we can say that \ref{pravg} is valid for $s=0$ as well.
\end{rem}

\begin{prp}$H$ is round regarding ${\cal{M}}^{acc}$ iff $l_1=l_2,\ |(H^{-k})^{(l_1)}|=|(H^{+k})^{(l_2)}|$ where $lev(H^{-k})=l_1, lev(H^{+k})=l_2, k={\cal{M}}^{acc}(H)$.
\end{prp}
\P Let $H^{(l)}\ne\emptyset, H^{(l+1)}=\emptyset,\A(H^{(l)})=k$. Obviously $(H^{-k})^{(l)},(H^{+k})^{(l)}$ are finite and non-empty i.e. $l=l_1=l_2$. Clearly $H$ is round iff
$$\frac{\A((H^{-k})^{(l)})+\A((H^{+k})^{(l)})}{2}=k$$
which is equivalent to $|(H^{-k})^{(l)}|=|(H^{+k})^{(l)}|$ by Proposition \ref{paround}.
\Pes

\begin{prp}$H$ is round regarding ${\cal{M}}^{lis}$ iff 
$$\frac{\varlimsup H^{-k}+\varliminf H^{+k}}{2}=\frac{\varlimsup H+\varliminf H}{2} \text{ where } k={\cal{M}}^{lis}(H).$$
\end{prp}
\P Let $a=\varliminf H=\varliminf H^{-k},d=\varlimsup H=\varlimsup H^{+k},b=\varlimsup H^{-k},c=\varliminf H^{+k}$. Roundness is equivalent to
$$\frac{\frac{a+b}{2}+\frac{c+d}{2}}{2}=\frac{a+d}{2}.$$
\Pes

\begin{prp}$H$ is round regarding ${\cal{M}}^{iso}$ iff 

1. ${\cal{M}}^{iso}(H^{-k})={\cal{M}}^{iso}(H^{+k})=k$ 

or

2. $\lim_{n\to\infty}\frac{|P_n|}{|S_n|}=1$

holds where $k={\cal{M}}^{iso}(H), S_n=\{x\in H-S(H',\frac{1}{n}):x\leq k\}, P_n=\{x\in H-S(H',\frac{1}{n}):x\geq k\}\ (n\in\mathbb{N}).$
\end{prp}
\P Let $s_n=\sum_{x\in S_n}x,p_n=\sum_{x\in P_n}x$. By definition $\lim_{n\to\infty}\A(S_n\cup P_n)=k.$ $H$ is round iff 
$$\lim_{n\to\infty}\frac{\A(S_n)+\A(P_n)}{2}=k.$$
That is equivalent to
$$\lim_{n\to\infty} \frac{\frac{s_n}{|S_n|} + \frac{p_n}{|P_n|}}{2} - \frac{s_n+p_n}{|S_n|+|P_n|}=0$$
(we make the remark that if $k\in S_n\cap P_n$ then the formula gets slightly different however in the limit it does not make any difference).

That is the same that
$$\frac{1}{2}\frac{s_n}{|S_n|}\frac{|P_n|-|S_n|}{|S_n|+|P_n|} + \frac{1}{2}\frac{p_n}{|P_n|}\frac{|S_n|-|P_n|}{|S_n|+|P_n|}\to 0$$
or
$$\frac{1}{2}\frac{|P_n|-|S_n|}{|S_n|+|P_n|}\Big(\frac{s_n}{|S_n|}-\frac{p_n}{|P_n|}\Big)\to 0.$$
Obviously 
$$-1<\frac{|P_n|-|S_n|}{|S_n|+|P_n|}=1-\frac{2}{1+\frac{|P_n|}{|S_n|}}<1$$
and 
$$\frac{s_n}{|S_n|}\leq k\leq \frac{p_n}{|P_n|}.$$
Because $\frac{s_n}{|S_n|},\frac{p_n}{|P_n|}$ both convergent we have two options:

1. They both converge to $k$ or

2. $\frac{|P_n|}{|S_n|}\to 1$.
\Pes

\begin{rem}For this mean ${\cal{M}}^{iso}$ in the definition of roundness it is important that $H^{-k},H^{+k}\in Dom({\cal{K}})$ as one can easily create an example where $\frac{s_n}{|S_n|},\frac{p_n}{|P_n|}$ are not convergent while $\frac{s_n+p_n}{|S_n|+|P_n|}$ do converge.
\end{rem}


{\footnotesize

\noindent

\noindent E-mail: alosonczi1@gmail.com\\

\end{document}